\documentclass[11pt]{article}

\usepackage[english]{style}
\usepackage{latexsym,amsmath,amssymb,amsfonts,graphics,amscd}
\usepackage[all]{xy}
\usepackage{pgf,tikz, tikz-cd}
\usetikzlibrary{arrows}
\usepackage{amssymb}
\DeclareFontFamily{U}{rsf}{}
\DeclareFontShape{U}{rsf}{m}{n}{
  <5> <6> rsfs5 <7> <8> <9> rsfs7 <10-> rsfs10}{}
\DeclareMathAlphabet{\mathscr}{U}{rsf}{m}{n}

\DeclareMathAlphabet{\mathgth}{U}{euf}{m}{n}

\DeclareFontFamily{U}{cyr}{}
\DeclareFontShape{U}{cyr}{m}{n}{
  <5> wncyr5 <6> wncyr6 <7> wncyr7 <8> wncyr8 <9> wncyr9 <10-> wncyr10}{}
\DeclareMathAlphabet{\mathcyr}{U}{cyr}{m}{n}

\input cyracc.def

\makeatletter
\def\operator@font{\sf}
\makeatother

\newcommand{\sbt}{{\scalebox{0.5}{\textbullet}}}
\newcommand{\eps}{\epsilon}

\newcommand{\sA}{{\mathcal A}}
\newcommand{\sDiff}{{\mathcal D iff}}

\newcommand{\sC}{{\mathcal C}}
\newcommand{\sD}{{\mathcal D}}

\newcommand{\cO}{{\mathscr O}}

\newcommand{\D}{{\mathbf D}}
\newcommand{\Diff}{{\mathbf Diff}}

\DeclareMathOperator{\Spec}{Spec}

\DeclareMathOperator{\Coh}{Coh}
\DeclareMathOperator{\Pol}{Pol}

\newcommand{\ra}{\rightarrow}

\newcommand{\C}{\mathbb{C}}

\newcommand{\iso}{\cong}

\newcommand{\del}{\partial}

\newcommand{\field}[1]{\mathbb{#1}}

\renewcommand{\L}{\field{L}}
\renewcommand{\phi}{\varphi}

\title{Weak $E_2$-Morita equivalences via quantization of the 1-shifted cotangent bundle}
\author{M\'arton Hablicsek\thanks{Leiden University, Netherlands
{\em e-mail: hablicsekhm@math.leidenuniv.nl}{\tt }}}
\begin{document}
\maketitle
In this paper, we investigate the structure of the convergent quantization of the 1-shifted cotangent bundle $S$ of a smooth scheme $X$ over a perfect field of positive characteristic. The quantization is an $E_2$-algebra over the Frobenius twist $S'$ of the 1-shifted cotangent bundle which restricted to the zero section $X'\ra S'$ is weakly $E_2$-Morita equivalent to the structure sheaf $\cO_{X'}$ of the Frobenius twist $X'$ of $X$. 

Explicitly, we show that the $(\infty,2)$-category of coherent (left-)modules over $\cO_{X'}$ is equivalent to the full subcategory of the $(\infty,2)$-category of coherent (left-)modules over the quantization restricted to the zero section are equivalent.
\section{Introduction}

\paragraph Deformation quantizations have interesting features over fields of positive characteristic. For instance, the Poisson-center of a Poisson-variety $(S,\{,\})$ is large: for any local sections $f$ and $g$ of $\cO_S$ we have $\{f^p,g\}=0$. Moreover, any convergent (central) quantization $\sA$ of $\cO_S$ is an Azumaya algebra (see \cite{BezKal} for more details). In the simplest case, when $S=T^*X$ is the cotangent bundle of a smooth variety $X$ equipped with its natural symplectic structure, the convergent quantization is the sheaf of rings of crystalline differential operators, $\D_X$. It is an Azumaya algebra over the Frobenius twist $S'=T^*X'$. Moreover, it is naturally a trivial Azumaya algebra restricted to the zero section $X'\rightarrow T^*X'=S'$ of the Frobenius twist of the cotangent bundle (see \cite{BezMirRum}) meaning that $\D_X|_{X'}$ is Morita equivalent ot $\cO_{X'}$. 

The primary goal of this paper is to generalize the above and formulate a conjecture (and provide some evidence) about the structure of deformation quantizations of shifted symplectic structures in positive characteristics. 

\paragraph Recently, a new approach to deformation quantization was introduced by Calaque, Pantev, To\"en, Vaqui\'e and Vezzosi (\cite{PTVV,CPTVV}) and separately by Pridham (\cite{Pri1,Pri2}) to study quantizations in the context of derived algebraic geometry. The input of the quantization problem is an $n$-shifted symplectic derived stack $S$, the output is a sheaf of $B\D_{n+1}$-algebras (deforming the structure sheaf of $S$). 

\paragraph Consider the simplest, non-trivial case, the $-1$-shifted cotangent bundle $T^*[-1]X$ of a smooth Calabi-Yau variety $X$ of dimension $\dim X<p$. The $-1$-shifted cotangent bundle is a derived scheme which can be obtained as the derived self-intersection of the zero section $X\ra T^*X$:
\[\xymatrix{T^*[-1]X\ar[r]\ar[d]& X\ar[d]^0 \\ X\ar[r]^{0} & T^*X}.\]
The structure sheaf, $\cO_{T^*[-1]X}$, (over $X$) is quasi-isomorphic to the derived tensor product
\[\cO_X\otimes^L_{\cO_{T^*X}} \cO_X.\]
Given a volume form $\omega$ on $X$, we have a trivialization $\omega_X\iso \cO_X$. Using this trivialization $\omega_X\iso \cO_X$, the complex $\cO_X\otimes^L_{\cO_{T^*X}} \cO_X$ is quasi-isomorphic (up to a shift) to the formal complex
\[0\ra \cO_X\xrightarrow{0}\Omega^1_X\xrightarrow{0} \Omega^2_X\ra...\]
The convergent quantization using the Costello-Li framework (\cite{CosLi}) is a Batalin-Vilkovisky algebra structure on the structure sheaf of the $T^*[-1]X$. In our case, the convergent quantization is the de Rham complex of $X$ (up to a shift) where the Batalin-Vilkovisky-differential is basically given by the de Rham differential (see \cite{BarGin} and \cite{BehFan} for more details). In positive characteristics, the de Rham complex is naturally a complex of $\cO_{X'}$-module, hence naturally a complex over the Frobenius twist $X'$. The Cartier isomorphisms identify the cohomology sheaves of the de Rham complex: the $d$-th cohomology sheaf is the sheaf of algebraic $d$-forms on $X'$. Moreover, the de Rham complex of $X$ can be realized as a line bundle (up to a shift) on the Frobenius twist, $T^*[-1]X'$ of the $-1$-shifted cotangent bundle, and this line bundle becomes a trivial line bundle once restricted to the zero section $X'\rightarrow T^*[-1]X'$ (see \cite{AriCalHab}, \cite{Hab2} for more details). 

We summarize the discussion above as follows.
\[\begin{array}{c|c|c}
\mbox{Shift }(n) & \mbox{Structure of the quantization}& \mbox{Structure over the zero section }\\
& \mbox{over the Frobenius twist} & X'\rightarrow T^*[n]X' \\ \hline
-1 & \mbox{Line bundle} & \mbox{Trivial}\\
0 & \mbox{Azumaya algebra} & \mbox{Trivial}
\end{array}\]

\paragraph From the table above we can see that the structure of the quantization gets more interesting as $n$ increases. It is natural to ask the following question.

\begin{center}\textbf{What kind of structure does the convergent quantization of an $n$-shifted symplectic derived stack have in positive characteristics?}
\end{center} 

Generalizing from the table above, the natural answer to this question would be that the quantization is an $E_{n+1}$-Azumaya algebra over the Frobenius twist, which becomes trivial (an $E_{n+1}$-Morita equivalence) once we restrict it to the zero section $X'\rightarrow T^*[n]X'$. The purpose of this paper is to investigate the structure of the convergent quantization of the 1-shifted cotangent bundle, $T^*[1]X$, restricted to the zero section $X'\rightarrow T^*[1]X'$.

\paragraph The notion of $E_n$-Morita equivalence of algebras was recently defined by Haugseng (\cite{Hau}) (or for the ``pointed'' version, see \cite{GwiSch}, \cite{Sch}). For nice enough $(\infty,1)$-category $\sC$ the $E_n$-algebras over $\sC$ form an $(\infty, n+1)$-category. This category can be described roughly as follows: objects are the $E_n$-algebras, $1$-morphisms are $E_{n-1}$ algebras in bimodules, $2$-morphisms are $E_{n-2}$-algebras in bimodules over the bimodules, etc. This notion recovers the standard Morita 2-category whose objects are associative algebras, 1-morphisms are bimodules and 2-morphisms are morphisms of bimodules (see \cite{AriCalHab} for Azumaya schemes and \cite{Toe} for derived schemes).

\paragraph Explicitly two $E_2$-algebras $A$ and $B$ are $E_2$-Morita equivalent if there exist algebras $N$ and $M$ in $A-B$ and $B-A$ bimodules so that $N\otimes^L_A M$ is $E_1$-Morita equivalent to $B$ and $M\otimes^L_B N$ is $E_1$-Morita equivalent to $A$. It is easy to see that if $A$ is commutative algebra, and $B$ is $E_2$-Morita equivalent to $A$, then $B$ has to be perfect as an $A$-module.

\paragraph In our set-up, we consider the $1$-shifted cotangent bundle $T^*[1]X$ of a smooth scheme $X$ over a perfect field of characteristic $p>\dim X$. The convergent quantization of $\cO_{T^*[1]X}$ is a variant of the Hochschild cosimplicial complex, which we call the crystalline Hochschild cosimplicial complex. Its restriction to the zero section $X'\ra T^*[1]X'$ is quasi-isomorphic to the $\cO_{X'}$-linear Hochschild cosimplicial complex of $\cO_X$ (Proposition \ref{prop:rest}), which we denote by $\sDiff_{\cO_{X'}}(\cO_X^\sbt,\cO_X)$.\footnote{Here we regard $\cO_X$ as a $\cO_{X'}$-algebra using the Frobenius morphism $F:X\ra X'$. We abuse notation and we will write $\cO_X$ instead of $F_*\cO_X$.} A quick consequence of the above paragraph is that $\cO_{X'}$ and $\sDiff_{\cO_{X'}}(\cO_X^\sbt,\cO_X)$ are NOT $E_2$-Morita equivalent as $\sDiff_{\cO_{X'}}(\cO_X^\sbt,\cO_X)$ is not perfect as a $\cO_{X'}$-complex.

\paragraph On the other hand, we show that the $(\infty,2)$-category of coherent (left-)modules over $\cO_{X'}$ is equivalent to a full subcategory of the $(\infty,2)$-category of coherent (left-)modules over $\sDiff_{\cO_{X'}}(\cO_X^\sbt,\cO_X)$ generated by $\cO_X$. Explicitly, in Proposition \ref{prop:quant2}, we show that the algebra $\cO_X$ can be realized as a brace module over $\sDiff_{\cO_{X'}}(\cO_X^\sbt,\cO_X)$. This construction provides functors $F$ and $G$ between the category of coherent (left-)modules over $\cO_{X'}$ and over $\sDiff_{\cO_{X'}}(\cO_X^\sbt,\cO_X)$ as follows:
\[F:\Coh(\sDiff_{\cO_{X'}}(\cO_X^\sbt,\cO_X))\ra \Coh(\cO_{X'})\quad F(-)=\cO_X\otimes^{\L}_{\sDiff_{\cO_{X'}}(\cO_X^\sbt,\cO_X)}-\]
\[G:\Coh(\cO_{X'})\ra \Coh(\sDiff_{\cO_{X'}}(\cO_X^\sbt,\cO_X))\quad G(-)=\cO_X\otimes^{L}_{\cO_{X'}}-.\]
We prove the following theorem relating the categories of coherent (left-)modules over $\cO_{X'}$ and $\sDiff_{\cO_{X'}}(\cO_X^\sbt,\cO_X)$.

\begin{Theorem}\label{thm:main1}
The following two $(\infty,2)$-categories
\begin{itemize}
\item the category of coherent (left-)modules over $\cO_{X'}$ where we consider $\cO_{X'}$ as an $E_2$-algebra and
\item the full thick subcategory of the coherent (left-)modules over $\sDiff_{\cO_{X'}}(\cO_X^\sbt,\cO_X)$ generated by $\cO_X$
\end{itemize}
are equivalent. 
\end{Theorem}
\medskip

Explicitly, we show that for any algebra object $A$ in the categories $\Coh(\cO_{X'})$ we have that
$FG(A)$ and $A$ are Morita equivalent.

\paragraph \textbf{Question:} Is the full thick subcategory of the category of coherent (left-)modules over $\sDiff_{\cO_{X'}}(\cO_X^\sbt,\cO_X)$ generated by $\cO_X$ actually equivalent to the $(\infty,2)$-category of coherent (left-)modules over $\sDiff_{\cO_{X'}}(\cO_X^\sbt,\cO_X)$?

We also expect that our theorem generalizes for higher shifts as well.
\begin{Conjecture}
Let $S$ be the $n$-shifted cotangent bundle of a smooth scheme $X$ over a perfect field $k$ of characteristic $p$. Let $\sA$ be the convergent quantization of $\cO_S$. Consider the Frobenius twist $S'$ of $S$, and the zero section $i:X'\ra S'$. Then, the algebra $\sA$ can be regarded as an $E_{n+1}$-algebra over $S'$ so that
\begin{itemize}
\item Weak Morita equivalence: The $(\infty, n+1)$-category of coherent sheaves over $\cO_{S'}$ (viewed as an $E_{n+1}$-algebra) is equivalent to a full thick subcategory of the $(\infty, n+1)$-category of coherent modules over $i^*\sA$; and 
\item Weak Azumaya property: \'Etale locally over $X$, the $(\infty, n+1)$-category of coherent sheaves over $\cO_{S'}$ (viewed as an $E_{n+1}$-algebra) is equivalent to a full thick subcategory of the $(\infty, n+1)$-category of coherent modules over $\sA$.
\end{itemize}
\end{Conjecture}

\paragraph We also remark that in the case of $-1$-shifted symplectic derived Artin stacks, the convergent quantization is a $p$-torsion element of the Picard group (for the de-Rham complex, see \cite{Hab2}), in the case of symplectic varieties, the convergent quantization is a $p$-torsion element of the Brauer group (\cite{BezKal}). We wonder whether the quantizations in higher shifts can be realized as $p$-torsion elements of higher Brauer groups (\cite{Hau}).

\paragraph The paper is organized as follows. In Section \ref{sec:two}, we collect facts about the Hochschild cosimplicial complex and its variant which we call the crystalline Hochschild cosimplicial complex. We show (Proposition \ref{prop:rest}) that the $\cO_{X'}$-linear Hochschild cosimplicial complex $\sDiff_{\cO_{X'}}(\cO_{X}^\sbt,\cO_X)$ is the restriction of the convergent quantization of $T^*[1]X$ to the zero section $X'\ra T^*[1]X'$. We also provide a left and right brace module structure of $\cO_X$ over the brace algebra $\sDiff_{\cO_{X'}}(\cO_{X}^\sbt,\cO_X)$ (Proposition \ref{prop:quant2}). In Section \ref{sec:main}, we prove our main theorem, Theorem \ref{thm:main1}.

\paragraph[Acknowledgment:] The author thanks Rune Haugseng for patiently explaining his work and answering numerous questions. The author also thanks Ryszard Nest, Pavel Safronov, and Elden Elmanto for various suggestions and fruitful discussions.
\section{The two Hochschild cosimplicial complexes}\label{sec:two}

In this section, we recollect facts about two cosimplicial complexes: the first one, $\sDiff(\cO_X^\sbt,\cO_X)$, is a global model for the Hochschild cosimplicial complex equipped with its brace algebra structure, the second one, $\Diff(\cO_X^\sbt,\cO_X)$, can be thought of the convergent quantization of the 1-shifted cotangent bundle $T^*[1]X$ equipped with its brace algebra structure. Similarly, for every $\cO_X$-bimodule, $P$, we define the complexes $\sDiff(\cO_X^\sbt,P)$ and $\Diff(\cO_X^\sbt, \cO_X)$, which are cosimplicial complexes equipped with a left brace module structure structure over the Hochschild complexes (see \cite{CalvdB}, \cite{CPTVV}, \cite{Hab}, \cite{Yek} for more details). Finally, we prove that the $\cO_{X'}$-linear Hochschild cosimplicial complex $\sDiff_{\cO_{X'}}(\cO_{X}^\sbt,\cO_X)$ is the restriction of $\Diff(\cO_X^\sbt, \cO_X)$ to the zero section $X'\ra T^*[1]X'$. We also provide a left and right brace module structure of $\cO_X$ over the brace algebra $\sDiff_{\cO_{X'}}(\cO_{X}^\sbt,\cO_X)$.

\subsection{Grothendieck polydifferential operators}
\label{sec:grpo}
We begin with the definition of Grothendieck (poly)differential operators.

\begin{Definition}
Let $P$ be an $\cO_X$-bimodule, and $A:\cO_X\ra P$ a $k$-linear map. Given a sequence of functions $f_0, f_1, ..., \in \cO_X$, define a sequence of $k$-linear maps $A_m:\cO_X\ra P$ given by $A_{-1}=A$ and, $A_n:= f_n A_{n-1}-A_{n-1}f_n$. We say that $A$ is a differential operator of order at most $N$ if for every point $x\in X$ and every section $s\in \cO_X$ defined at $x$, there exists a neighborhood $U$ of $x$ and $N\geq 0$ such that for any open subset $V\subset U$ and any choice of functions $f_0, f_1,..., f_N$ on $V$, so that $A_N(s|_V)$ vanishes.
\end{Definition}
\medskip 

The Grothendieck differential operators $\cO_X \ra  P$ form a sheaf that we denote by $\sDiff_X(\cO_X, P)$. In the case of $P=\cO_X$ we denote the sheaf of differential operators by $\sD_X:=\sDiff_X(\cO_X,\cO_X)$. The tangent bundle $T_X$ of $X$ has a $(k,\cO_X)$ Lie-algebroid structure. The ring $\sD_X$ is the universal PD-enveloping algebra of the $(k,\cO_X)$ Lie-algebroid $T_X$. For instance, if $X=\Spec k[x]$, then $\sD_X$ is the PD-polynomial ring of one variable, $k\langle x \rangle$. Moreover, by definition the ring $\sD_X$ comes with a filtration
\[\cO_X=\sD_X^{\leq 0}\subset \sD_X^{\leq 1}\subset \sD_X^{\leq 2}\subset...\]
given by the degree of the differential operators.

\begin{Definition} A $k$-polylinear map $A: \cO_X\times...\times \cO_X\ra P$ (of $n$ arguments) is a polydifferential
operator of (poly)order at most $(N_1,..., N_n)$ if it is a differential operator of order at most $N_j$ in the $j$-th argument whenever the remaining $n-1$ arguments are fixed.
\end{Definition}
\medskip

The polydifferential operators $\cO_X\times...\times \cO_X\ra P$ of $i$ arguments form a sheaf, which we denote by $\sDiff(\cO_X^i,P)$.

\paragraph We can identify the sheaf $\sDiff(\cO_X^i,P)$ with the tensor product
\[\sD_X\otimes_{\cO_X}....\otimes_{\cO_X}\sD_X\otimes_{\cO_X}P\]
as follows (here the number of the $\sD_X$ terms is $i$). The map
\[\sD_X\otimes_{\cO_X}....\otimes_{\cO_X}\sD_X\otimes_{\cO_X}P\ra \sDiff(\cO_X^i,P)\]
given by
\[A_1\otimes...\otimes A_i\otimes p\mapsto A_1(-)A_2(-)...A_i(-)p\]
(for local sections $A_i\in \sD_X$ and $p\in P$) is clearly an isomorphism. (Here we use the natural $\cO_X$-bimodule structure on $\sD_X$.)

\subsection{The Grothendieck Hochschild cosimplical complex}
\label{sec:ghc}

\paragraph The sheaves of polydifferential operators form a natural cosimplicial complex $\sDiff(\cO_X^\sbt,P)$ whose $i$-th term is $\sDiff(\cO_X^i,P)$ and the differentials $d_{i,j}$ are given by 
\begin{align*}
d_{i,k}A(g_1,...,g_{i+1})=
\begin{cases}
g_1A(g_2,...,g_{i+1}) & k=0\\
A(g_1,...,g_kg_{k+1},...,g_{i+1}) & 0<k<i+1\\
A(g_1,...,g_i)g_{i+1} & k=i+1
\end{cases}
\end{align*}

where $A:\cO_X\times...\times \cO_X\ra P$ is a polydifferential operator of $i$ arguments and $\{g_1,...,g_{i+1}\}$ is a local section of $\cO_X\times...\times \cO_X$ (of $i+1$ arguments).

\paragraph Let $A\in \sDiff(\cO_X^i,\cO_X)$ and $B\in \sDiff(\cO_X^j,\cO_X)$. We define $A\cdot B\in \sDiff(\cO_X^{i+j},\cO_X)$ as the differential operator mapping $a_1,...,a_{i+j}$ to 
\[(-1)^{ij}A(a_1,...,a_i)\cdot B(a_{i+1},...,a_{i+j}).\]
This product endows the cosimplicial complex $\sDiff(\cO_X^\sbt,\cO_X)$ a cosimplicial algebra structure.

\paragraph Similarly, we see that if $P$ is an $\cO_X$-bimodule, then the complex $\sDiff(\cO_X^\sbt,P)$ has a bimodule structure over $\sDiff(\cO_X^\sbt,\cO_X)$. Moreover, if the $\cO_X$-bimodule, $P$ has an associative algebra structure, then $\sDiff(\cO_X^\sbt,P)$ has a cosimplicial algebra structure defined parallel to the cosimplicial algebra structure on $\sDiff(\cO_X^\sbt,\cO_X)$.

\paragraph The cosimplicial algebra, $\sDiff(\cO_X^\sbt,\cO_X)$ is equipped with a brace algebra structure as follows. Let $A\in \sDiff(\cO_X^i,\cO_X)$ and $A_{l}\in \sDiff(\cO_X^{j_l},\cO_X)$ (for $l=1,...,m$). The brace operations $A\{A_1,...,A_m\}$ are defined as operations of degree $-m$, i.e. $A\{A_1,...,A_m\}\in \sDiff(\cO_X^n,\cO_X)$, where $n=i+\sum_{l=1}^m j_l-m$. Explicitely, $A\{A_1,...,A_m\}$ maps $n$ sections of $\cO_X$, $a_1,...,a_n$ to
\[\sum_{0\leq i_1\leq...\leq i_m\leq n}(-1)^\epsilon A(a_1,...,a_{i_1},A_1(a_{i_1+1},...),...,a_{i_m},A_m(a_{i_m+1},...),...,a_n)\]
where $\eps:=\sum_{i=1}^m i_l(j_l-1)$. 

Similarly, if $P$ is an $\cO_X$-bimodule, then $\sDiff(\cO_X^\sbt,P)$ has a (left) module structure over the brace algebra $\sDiff(\cO_X^\sbt,\cO_X)$ as follows. Let $B\in \sDiff(\cO_X^\sbt,P)$ and $A_{l}\in \sDiff(\cO_X^{j_l},\cO_X)$ (for $l=1,...,m$). The brace operations $B\{A_1,...,A_m\}$ are defined as operations of degree $-m$, i.e. $B\{A_1,...,A_m\}\in \sDiff(\cO_X^n,P)$, where $n=i+\sum_{l=1}^m j_l-m$. Explicitly, $B\{A_1,...,A_m\}$ maps $n$ sections of $\cO_X$, $a_1,...,a_n$ to
\[\sum_{0\leq i_1\leq...\leq i_m\leq n}(-1)^\epsilon B(a_1,...,a_{i_1},A_1(a_{i_1+1},...),...,a_{i_m},A_m(a_{i_m+1},...),...,a_n)\]
where $\eps:=\sum_{i=1}^m i_l(j_l-1)$. 

\paragraph Similarly, a right $\sDiff(\cO_X^\sbt,\cO_X)$ brace module structure on a complex is equivalent to a left $\sDiff(\cO_X^\sbt,\cO_X)^{op}$ brace module structure. 

\paragraph Note that the Grothendieck ring of differential operators $\sD_X$ is in general not a finitely generated algebra over $\cO_X$. On the other hand, the homotopy groups of the Grothendieck complex are locally free $\cO_X$-modules by the Hochschild-Konstant-Rosenberg theorem (\cite{Swa, Yek}) if the characteristic $p$ is large enough ($p\geq \dim X$).

\subsection{The crystalline differential operators}

\begin{Definition} The crystalline ring of differential operators, $\D_X=\Diff(\cO_X,\cO_X)$ is defined as the universal enveloping $(k,\cO_X)$ Lie-algebroid of the tangent bundle $T_X$ . 
\end{Definition}
\medskip

\paragraph \textbf{Example:} in the case of $X=\Spec k[x]$, the ring of differential operators $\D_X$ is the Weyl-algebra 
\[k\langle x,\frac{d}{dx}\rangle/(\frac{d}{dx}x-x\frac{d}{dx}-1).\]

\paragraph The algebra $\D_X$ is equipped with a filtration 
\[\cO_X=\D_X^{\leq 0}\subset \D_X^{\leq 1}\subset \D_X^{\leq 2}\subset...\]
given by the degree of the differential operator. Note that the filtered pieces $\D_X^{\leq p-1}$ and $\sD_X^{\leq p-1}$ are isomorphic as $\cO_X$-bimodules! Given a derivation $t\in T_X$, its $p$-th composite $t^{[p]}:=t\circ t\circ...\circ t$ is also a derivation. This gives rise to a distinguished element $t^p-t^{[p]}$ in $\D_X$ for every derivation $t\in T_X$. The quotient of $\D_X$ with the ideal generated by these distinguished elements can be identified with $\D_X^{\leq p-1}$. Hence $\D_X^{\leq p-1}$ is a split sub-$\cO_X$-bimodule of $D_X$. This also provides an algebra structure on $\D_X^{\leq p-1}$ induced by the algebra structure on $\D_X$. Moreover, we have that
 
 \begin{Proposition}\cite{BezMirRum}
 The algebra $\D_X^{\leq p-1}$ is Morita equivalent to the structure sheaf $\cO_{X'}$. 
 \end{Proposition}

\paragraph Similarly as in Section \ref{sec:grpo}, given an $\cO_X$-bimodule $P$, we define the crystalline differential operators $\Diff(\cO_X,P)$ as the sheaf $\D_X\otimes_{\cO_X} P$ equipped with its natural bimodule structure. Moreover, we define the the polydifferential operators $\Diff(\cO_X^i,P)$ as the tensor product
\[\D_X\otimes_{\cO_X}....\otimes_{\cO_X}\D_X\otimes_{\cO_X}P\]
(here the number of the $\D_X$ terms is $i$).

\subsection{The crystalline Hochschild cosimplical complex}

\paragraph Similarly as in Section \ref{sec:ghc} we can construct a cosimplical complex $\Diff(\cO_X^\sbt,P)$ from the sheaves $\Diff(\cO_X^i,P)$. We call this cosimplical complex the \textit{crystalline Hochschild complex} of $P$. 

\paragraph[Remark:] Even though the crystalline ring of differential operators is a finitely generated algebra over $\cO_X$, the homotopy groups of the crystalline Hochschild cosimplicial complex are not locally free sheaves over $\cO_X$. 
\\

Hence, the two Hochschild cosimplicial complexes have peculiar behaviour. One consists of more complicated terms (terms obtained from non-finitely generated algebras) with less complicated homotopy groups (locally free sheaves), the other consists of less complicated terms (terms obtained from finitely generated algebras) with more complicated homotopy groups.

\paragraph The cosimplicial complex $\Diff(\cO_X^\sbt,\cO_X)$ can be thought of as the convergent quantization of the 1-shifted cotangent bundle $T^*[1]X$. We represent the 1-shifted cotangent bundle as $B(T^*X/X)$, i.e as the simplicial scheme
\[\begin{tikzcd}
...T^*X\times_X T^*X  \arrow[r]
\arrow[r, shift left=2]
\arrow[r, shift right=2] & T^*X \arrow[r, shift left]
\arrow[r, shift right]& X
\end{tikzcd}\]
where locally the maps are given by 
\begin{align*}
d_{i,k}(x,\omega_1,...,\omega_{i})=
\begin{cases}
(x,\omega_2,...,\omega_{i}) & k=0\\
(x,\omega_1,...,\omega_k+\omega_{k+1},..., \omega_{i}) & 0<k<i\\
(x,\omega_1,...,\omega_{i-1}) & k=i
\end{cases}
\end{align*}
(here we think of $T^*X$ as a group scheme over the base scheme $X$, where the fiber-wise group structure is given by the Abelian group structure on $\Omega^1_X$).  

Hence, the structure sheaf $\cO_{T^*[1]X}$ of $T^*[1]X$ is the represented by the cosimpicial algebra
\[\begin{tikzcd}
\cO_X \arrow[r, shift left]
\arrow[r, shift right] & \cO_{T^*X} \arrow[r]
\arrow[r, shift left=2]
\arrow[r, shift right=2] & \cO_{T^*X}\otimes_{\cO_X} \cO_{T^*X}...
\end{tikzcd}\]
which we will denote by $\Pol^\sbt(X)$ where the differentials are provided by the maps representing $T^*[1]X$ as the simplicial scheme above. We denote the Frobenius twist of the cosimplicial algebra $\Pol^\sbt(X)$ by $\Pol^\sbt(X')$, it is the cosimplicial algebra given by the Frobenius twists of the algebras above
\[\begin{tikzcd}
\Pol^\sbt(X')=\cO_{X'} \arrow[r, shift left]
\arrow[r, shift right] & \cO_{T^*X'} \arrow[r]
\arrow[r, shift left=2]
\arrow[r, shift right=2] & \cO_{T^*X'}\otimes_{\cO_{X'}} \cO_{T^*X'}...
\end{tikzcd}.\]

An important observation is that since the center of $\D_X$ can be identified with $\cO_{T^*X'}$ we have that $\Pol^\sbt(X')$ is a sub-complex of $\Diff(\cO_X^\sbt,\cO_X)$. Moreover, the brace structures on $\Diff(\cO_X^\sbt,\cO_X)$ become trivial on $\Pol^\sbt(X')$.

\paragraph The zero section $X\ra T^*[1]X$ is weak equivalent to the simplicial scheme
\[\begin{tikzcd}
...T^*X\times_X T^*X\times_X T^*X  \arrow[r]
\arrow[r, shift left=2]
\arrow[r, shift right=2] & T^*X\times_X T^*X \arrow[r, shift left]
\arrow[r, shift right]& T^*X
\end{tikzcd}\]
where locally the maps are given by
\begin{align*}d_{i,k}(x, \omega_1, ..., \omega_{i+1})=\begin{cases}
(x,\omega_2,...,\omega_{i}) & k=0\\
(x,\omega_1,...,\omega_k+\omega_{k+1},..., \omega_{i+1}) & 0<k<i+1\\
\end{cases},\end{align*}
(in other words, the scheme $X$ is weak-equivalent to $E(T^*X/X)$, the quotient of $T^*X$ with the action given by the group $T^*X$ over $X$ where the action is given by fiber-wise addition). 

The zero section $X\ra T^*[1]X$ can be equipped with a natural Lagrangian structure in the case when the base field is the field of complex numbers (see \cite{Cal} for more details). The following proposition can be thought of the quantization of this Lagrangian structure (\cite{Hab}). For the proof, we refer the reader to \cite{Sri}.

\begin{Proposition}\label{prop:quant} We have
\[\Diff(\cO_X^\sbt,\D_X)=\cO_X\]
(and similarly
\[\Diff(\cO_X^\sbt,\D_X^{op})^{op}=\cO_X).\]
\end{Proposition}

\subsection{Comparison results}

In this section, we present results comparing the Grothendieck and the Crystalline Hochschild complexes.

\paragraph Since the bracket on the tangent bundle $T_X$ is $\cO_{X'}$-linear, we can consider the universal enveloping $(\cO_{X'},\cO_X)$ Lie-algebroid $\sDiff_{\cO_{X'}}(\cO_X,\cO_X)$ and the corresponding Hochschild complexes $\sDiff_{\cO_{X'}} (\cO_X^\sbt,P)$. 

\begin{Lemma}\label{lem:sdiffperd}
We have an isomorphism of algebras
\[\sDiff_{\cO_{X'}}(\cO_X,\cO_X)\iso \D_X^{\leq p-1}.\]
\end{Lemma}

\begin{Proof}
For any non-trivial derivation $t\in T_X$, its $p$-th divided power is not $\cO_{X'}$-linear. Moreover for any derivation $t\in T_X$, $t^p$ and $t^{[p]}$ act indentically on local functions of $\cO_X$. As a consequence, $\sDiff_{\cO_{X'}}(\cO_X,\cO_X)$ is the quotient of $\D_X$ by the ideal generated by the $t^p-t^{[p]}$. \qed
\end{Proof}
\medskip

Now, we compare $\sDiff_{\cO_{X'}}(\cO_X^\sbt,\cO_X)$ with the crystalline Hochschild complex $\Diff(\cO_X^\sbt,\cO_X)$.

\begin{Proposition}\label{prop:rest}
The restriction of the convergent quantization of $T^*[1]X$ to the zero section $X'\ra T^*[1]X'$ is the $\cO_{X'}$-linear Hochschild cosimplicial complex of $\cO_X$:
\[\sDiff_{\cO_{X'}}(\cO_X^\sbt,\cO_X)=\Diff(\cO_X^\sbt,\cO_X)\otimes_{\Pol^\sbt(X')} \cO_{X'}.\]
\end{Proposition}

\begin{Proof}
First of all, $\Diff(\cO_{X}^i, \cO_X)$ is a free-module over the algebra $\Pol^i(\cO_{X'})$ locally generated by the polydifferential operators $\partial_1^{\alpha_1}\otimes \partial_2^{\alpha_2}\otimes...\otimes \partial_i^{\alpha_i}$ of $\D_X\otimes \D_X\otimes...\otimes \D_X$ with orders $\alpha_i<p$. Similarly, using Lemma \ref{lem:sdiffperd}, $\sDiff_{\cO_{X'}}(\cO_X^i,\cO_X)$ can be identified with the tensor product
\[\D_X^{\leq p-1}\otimes_{\cO_X}\D_X^{\leq p-1}\otimes_{\cO_X}...\otimes_{\cO_X}\D_X^{\leq p-1}.\]

Similarly, the complex $\sDiff_{\cO_{X'}}(\cO_{X}^\sbt, \cO_X)$ can be realized as a subcomplex of $\Diff(\cO_X^i,\cO_X)$ (but not a subalgebra) since the differential of $\Diff(\cO_X^i,\cO_X)$ respects the filtration. This subcomplex is the 
 $<p$ filtrered part of $\Diff(\cO_X^\sbt, \cO_X)$, identifying it with $\Diff(\cO_X^\sbt,\cO_X)\otimes_{\Pol^\sbt(X')} \cO_{X'}$.

This also shows that both the brace structures and differentials align.\qed
\end{Proof}
\medskip

Similar to Proposition \ref{prop:quant} we can equip the algebra $\cO_X$ with a left-(and right-)brace module structure over $\sDiff_{\cO_{X'}}(\cO_X^\sbt,\cO_X)$. 

\begin{Proposition}\label{prop:quant2} 
We have quasi-isomorphisms of complexes
\[\sDiff_{\cO_{X'}}(\cO_X^\sbt, \D^{\leq p-1}_X)=\cO_X\]
(and 
\[\sDiff_{\cO_{X'}}(\cO_X^\sbt, (\D^{\leq p-1}_X)^{op})^{op}=\cO_X).\]
\end{Proposition}

\begin{Proof}
The zeroth homotopy sheaf can be computed easily, it is the centralizer of $\cO_X$ inside $\D_X^{\leq p-1}$ which is $\cO_X$. In order to show that the higher homotopy sheaves vanish, we consider the map of cosimplicial complexes
\[\phi: \Diff(\cO_X^\sbt,\D_X)\ra \sDiff_{\cO_{X'}}(\cO_X^\sbt, \D^{\leq p-1}_X)\]
given by the identifications
\[\Diff(\cO_X^i,\D_X)=\D_X\otimes_{\cO_X}\D_X\otimes_{\cO_X}...\otimes_{\cO_X}\D_X,\]
\[\sDiff_{\cO_{X'}}(\cO_X^i, \D^{\leq p-1}_X)=\D_X^{\leq p-1}\otimes_{\cO_X}\D_X^{\leq p-1}\otimes_{\cO_X}...\otimes_{\cO_X}\D_X^{\leq p-1}\]
(with $i+1$ terms!) and the quotient map $\D_X\ra \D^{\leq p-1}_X$ given by the ideal generated by the distinguished elements $t^p-t^{[p]}$. More precisely, a crystalline polydifferential operator $\Diff(\cO_X^i, \D_X)$ can be considered as a $\cO_{X'}$-linear polydifferential operator providing a map $\Diff(\cO_X^i, \D_X)\ra \sDiff_{\cO_{X'}}(\cO_X^i, \D_X)$. This map is not injective, since $t^p$ and $t^{[p]}$ (with $t\in T_X$) act the same way on $\cO_X$. Finally, we compose this map with the quotient map $\sDiff_{\cO_{X'}}(\cO_X^i, \D_X)\ra \sDiff_{\cO_{X'}}(\cO_X^i, \D_X^{\leq p-1})$.

This map, $\phi$ splits, the inclusion map $i:\D_X^{\leq p-1}\ra \D_X$ gives an inverse. Moreover, the maps $i$ and $\phi$ respect the differentials of the Hochschild cosimplicial complexes, hence the complex $\sDiff_{\cO_{X'}}(\cO_X^\sbt, \D_X^{\leq p-1})$ is a direct summand of the complex $\Diff(\cO_X^\sbt,\D_X)$. The latter is quasi-isomorphic to $\cO_X$ by Proposition \ref{prop:quant}. Since the higher homotopy sheaves of $\Diff(\cO_X^\sbt,\D_X)$ vanish, the higher homotopy sheaves of $\sDiff_{\cO_{X'}}(\cO_X^\sbt, \D^{\leq p-1}_X)$ have to vanish as well.\qed
\end{Proof}
\medskip

Finally, we compute the homotopy sheaves of the Hochschild complex of $\D^{\leq p-1}_X$ over $\cO_X$ using the natural inclusion $\cO_X\ra \D^{\leq p-1}_X$ (which equips $\D^{\leq p-1}_X$ with a bimodule structure over $\cO_X$).

\begin{Proposition}\label{prop:dhoch}
\[\sDiff_{\cO_X}((\D^{\leq p-1}_{X})^\sbt,\D^{\leq p-1}_{X})=\cO_{X'}\]
\end{Proposition}

\begin{Proof}
First, we show that the zeroth homotopy sheaf is $\cO_{X'}$, and then we show that the higher homotopy sheaves vanish. The zeroth term of the complex is $\sDiff_{\cO_X}(\cO_X,\D^{\leq p-1}_X)$ which can be identified with $\cO_X$. Hence, the zeroth homotopy sheaf is the centralizer of $\D^{\leq p-1}_X$ inside $\cO_X$ which is indeed $\cO_{X'}$. 

Next, we show that the higher homotopy sheaves vanish. First, note that we only need to solve the problem locally, so we can assume that $X$ is a spectrum of a polynomial ring. Moreover, using K\"{u}nneth-formula, we can assume that the polynomial ring is of one variable.

We use the Dold-Kan correspondance, and we resolve $\D_X^{\leq p-1}=k\langle x, d\rangle/(dx-xd-1,d^p)$ with locally free $\D^e=\D_X^{\leq p-1}\otimes_{\cO_X} \D_X^{\leq p-1}$-modules. We claim that there is a 2-periodic resolution given by 
\[...\D^e\xrightarrow{d\otimes 1-1\otimes d}\D^e\xrightarrow{d^p-1\otimes 1+d^{p-2}\otimes d+...+1\otimes d^{p-1}}\D^e\xrightarrow{d\otimes 1-1\otimes d}\D^e\xrightarrow{m} \D_X^{\leq p-1}.\]
Here the first map is given by multiplication. The kernel of that map is the (left)ideal generated by $d\otimes 1-1\otimes d$ and by those monomials $d^i\otimes d^j$ for which $i+j\geq p$. However, the latter is also generated by the former. 

We turn our attention to the kernel of the map $\D^e\xrightarrow{d\otimes 1-1\otimes d}\D^e$ given by multiplication from the right by $d\otimes 1-1\otimes d$. Any element of $\D^e$ can be written as $\sum_{0\leq i,j\leq p-1} f_{i,j}d^i\otimes d^j$. We say that the degree of a monomial $f_{i,j}d^i\otimes d^j$ is $i+j$. 

The map $\D^e\xrightarrow{d\otimes 1-1\otimes d}\D^e$ sends the element $\sum_{0\leq i,j\leq p-1} f_{i,j}d^i\otimes d^j$ to the element
\[\sum_{0\leq i,j\leq p-1} (f_{i,j}d^{i+1}\otimes d^j-f_{i,j}d^i\otimes d^{j+1}).\]
Notice that the degree of each monomial is increased by 1. As a consequence, it is clear that if the element $\sum_{0\leq i,j\leq p-1} f_{i,j}d^i\otimes d^j$ is in the kernel of $\D^e\xrightarrow{d\otimes 1-1\otimes d}\D^e$, then degreewise it is in the kernel, i.e
\[\sum_{\substack{0\leq i,j\leq p-1\\ i+j=k}} f_{i,j}d^i\otimes d^j\]
is in the kernel as well for every $k$. Then, a simple calculation implies that $f_{i,j}=0$ if $i+j<p-1$, and that the kernel of $\D^e\xrightarrow{d\otimes 1-1\otimes d}\D^e$ is the (left)ideal generated by $d^p-1\otimes 1+d^{p-2}\otimes d+...+1\otimes d^{p-1}$. A similar calculation can be done for the kernel of $\D^e\xrightarrow{d^p-1\otimes 1+d^{p-2}\otimes d+...+1\otimes d^{p-1}}\D^e$.

We turn our attention now to compute 
\[\sDiff_{\cO_X}((\D^{\leq p-1}_{X})^\sbt,\D^{\leq p-1}_{X}).\]
We note that $\sDiff_{D^e}(\D^e, \D^{\leq p-1}_X)=k[x]$, because the image $t$ of $1\in \D^e$ has to satisfy $xt-tx=0$. As a consequence using our two-periodic resoluion the complex $\sDiff_{\cO_X}((\D^{\leq p-1}_{X})^\sbt,\D^{\leq p-1}_{X})$ becomes
\[k[x]\xrightarrow{d/dx} k[x]\xrightarrow{d^{p-1}/dx^{p-1}} k[x]\xrightarrow{d/dx}...\]
where the first map is derivation, and the second one is the $p-1$-st composite of the derivation. This complex is exact except at degree 0 concluding our proof.\qed
\end{Proof}
\medskip

We conclude this section by showing that the $\sDiff_{\cO_{X'}}(\cO_X^\sbt,\cO_X)$-linear Hochschild cosimplicial complex of $\cO_X$ can be identified with $\cO_{X'}$ with its trivial brace algebra structure.

\begin{Theorem}
We have a quasi-isomorphism of $E_2$-algebras
\[ \sDiff_{\sDiff_{\cO_{X'}}(\cO_X^\sbt,\cO_X)}(\cO_X^\sbt,\cO_X)=\cO_{X'}.\]
\end{Theorem}

\begin{Proof}
First, we use Proposition \ref{prop:quant2} to equip $\cO_X$ with a $\sDiff_{\cO_{X'}}(\cO_X^\sbt,\cO_X)$-structure. Then, the left hand-side of the Theorem becomes
\[\sDiff_{\sDiff_{\cO_{X'}}(\cO_X^\sbt,\cO_X)}(\sDiff_{\cO_{X'}}(\cO_X, \D_X^{\leq p-1})^\sbt,\sDiff_{\cO_{X'}}(\cO_X, \D_X^{\leq p-1})).\]
We can see that the $E_2$-module $\sDiff_{\cO_{X'}}(\cO_X, \D^{\leq p-1}_{X})$ is generated by $\D^{\leq p-1}_{X}$ over $\sDiff_{\cO_{X'}}(\cO_X^\sbt,\cO_X)$. Therefore the complex
\[\sDiff_{\sDiff_{\cO_{X'}}(\cO_X^\sbt,\cO_X)}(\sDiff_{\cO_{X'}}(\cO_X, \D^{\leq p-1}_{X})^i,\sDiff_{\cO_{X'}}(\cO_X, \D^{\leq p-1}_{X}))\]
is determined by the complex
\[\sDiff_{\cO_{X'}}((\D^{\leq p-1}_{X})^i,\sDiff_{\cO_{X'}}(\cO_X, \D^{\leq p-1}_{X})).\]
However, it is not a quasi-isomorphism, elements of the above complex do not give rise to elements of the complex
\[\sDiff_{\sDiff_{\cO_{X'}}(\cO_X^\sbt,\cO_X)}(\sDiff_{\cO_{X'}}(\cO_X, \D^{\leq p-1}_{X})^i,\sDiff_{\cO_{X'}}(\cO_X, \D^{\leq p-1}_{X})).\]
The requirement that the map has to be $\sDiff_{\cO_{X'}}(\cO_X^\sbt,\cO_X)$-linear shows that only elements of
\[\sDiff_{\cO_X}((\D^{\leq p-1}_{X})^i,\D^{\leq p-1}_{X}))\]
will give rise elements of the original complex, and as a consequence, it is quasi-isomorphic to the original complex. Finally, Proposition \ref{prop:dhoch} implies the statement of the Theorem.\qed
\end{Proof}
\section{Proof of Theorem \ref{thm:main1}}\label{sec:main}

In this section, we prove our main theorem, Theorem \ref{thm:main1}. We begin with the following statement.

\begin{Theorem} We have
\[\cO_X \otimes^{\L}_{\sDiff_{\cO_{X'}}(\cO_X^\sbt,\cO_X)} \cO_X=\D_X^{\leq p-1}\]
as algeras!
\end{Theorem}

\paragraph \label{par:rem} Before we prove our main theorem  we recall how the associative algebra structure is defined on
\[\cO_X \otimes^{\L}_{\sDiff_{\cO_{X'}}(\cO_X^\sbt,\cO_X)} \cO_X.\]
We follow (\cite{Saf}). Let $A$ be a brace algebra and $N$ (and $M$ resp.) be a right (and left resp.) brace module over $A$. The bar complex $T_\sbt(A[1])$ has a structure of a dg bialgebra (\cite{GerVor} ). Tthe one-sided bar complexes $N\otimes T_\sbt(A[1])$ and $T_\sbt(A[1])\otimes M$ carry a natural structure of a dg algebra compatible with the $T_\sbt(A[1])$ comodule structure, providing a dg-algebra stucture on the cotensor product 
\[N\otimes T\sbt(A[1]) \otimes^{T_\sbt(A[1])}T_\sbt(A[1])\otimes M\]
which is quasi-isomorphic to the derived tensor product $N\otimes_A^L M$. We remark that this dg algebra structure is compatible with the dg-algebra structure on $N\otimes_\C M$ meaning that the natural map 
\[N\otimes_\C M\ra N\otimes_A^\L M\]
is a dg-algebra map (\cite{Saf}). (Since $N$ and $M$ are right and left brace modules, they have an associative algebra structure.)

In the case of $\cO_X$ over $\sDiff_{\cO_{X'}}(\cO_X^\sbt,\cO_X)$, we constructed the right and left brace module structure in Proposition \ref{prop:quant2}. We are ready to prove our main theorem.

\begin{Proof}
This statement is very similar to Theorem 4.3 in \cite{Hab}. We highlight the key steps.

By Proposition \ref{prop:quant2}, we can equip $\cO_X$ with a left and right brace module structure over $\sDiff_{\cO_{X'}}(\cO_X^\sbt,\cO_X)$. Hence the complex
\[\cO_X \otimes^\L_{\sDiff_{\cO_{X'}}(\cO_X^\sbt,\cO_X)} \cO_X\]
can be represented as
\[\sD^{\sbt}(X):=\sDiff_{\cO_{X'}}(\cO_X^\sbt,(\D_X^{\leq p-1})^{op})^{op} \otimes^{\L}_{\sDiff_{\cO_{X'}}(\cO_X,\cO_X)}\sDiff_{\cO_{X'}}(\cO_X^\sbt,\D_X^{\leq p-1}).\]
Using the cup-product, we can identify the above complex with
\[\sDiff_{\cO_{X'}}(\cO_X^\sbt,\D_X^{\leq p-1}\otimes_{\cO_X}\D_X^{\leq p-1}).\]
Consider the map
\[ \phi:\D_X^{\leq p-1}\ra \D_X^{\leq p-1}\otimes_{\cO_X}\D_X^{\leq p-1}\]
given by 
\[\del \mapsto d(\del)+i(\del)\]
where $d$ is the Hochschild differential
\[\D_X^{\leq p-1}\ra \sDiff_{\cO_{X'}}(\cO_X,\D_X^{\leq p-1})= \D_X^{\leq p-1}\otimes_{\cO_X}\D_X^{\leq p-1}\]
and $i$ is the map 
\[\D_X^{\leq p-1}\ra \D_X^{\leq p-1}\otimes_{\cO_X}\D_X^{\leq p-1}\]
sending $\del$ to $1\otimes \del$. In other words, with the identification
\[\D_X^{\leq p-1}\otimes_{\cO_X}\D_X^{\leq p-1}\ra \sDiff_{\cO_{X'}}(\cO_X, \D_X^{\leq p-1})\]
the map $\phi$ can be realized as multiplication from the left
\[\phi(\del)=(x\mapsto \del x).\]

As a result, we obtain a map of complexes
\[\chi: \D_X^{\leq p-1}\ra \sDiff_{\cO_{X'}}(\cO_X^\sbt,\D_X^{\leq p-1}\otimes_{\cO_X}\D_X^{\leq p-1})=\sDiff_{\cO_{X'}}(\cO_X^\sbt, \sDiff_{\cO_{X'}}(\cO_X, \D_X^{\leq p-1}).\]
Indeed, the image of $\chi$ is in the kernel of the differential
\[\sDiff_{\cO_{X'}}(\cO_X, \D_X^{\leq p-1})\ra \sDiff_{\cO_{X'}}(\cO_X,\sDiff_{\cO_{X'}}(\cO_X, \D_X^{\leq p-1}))\]
since an element $f\in \sDiff_{\cO_{X'}}(\cO_X, \D_X^{\leq p-1})$ is in the kernel if and only if for every local section $x\in \cO_X(U)$ we have that $fx=xf$. This means that for any local section $y\in \cO_X(U)$ we have that $f(y)x=f(yx)$ meaning that $f$ has to be given by left multiplication by an element $\del \in \D_X^{\leq p-1}$.

Now, we show that the map $\chi$ is an algebra map. We will prove this by showing that \'etale locally the algebra structure on $\D_X^{\leq p-1}$ can be lifted to $\D_X^{\leq p-1}\otimes_{k}\D_X^{\leq p-1})$. This is sufficient, since the algebra structure on $\sDiff_{\cO_{X'}}(\cO_X^\sbt,\D_X^{\leq p-1}\otimes_{\cO_X}\D_X^{\leq p-1})$ is induced by the algebra structure on 
\[\sDiff_{\cO_{X'}}(\cO_X^\sbt,(\D_X^{\leq p-1})^{op})^{op} \otimes^{\L}_{k}\sDiff_{\cO_{X'}}(\cO_X^\sbt,\D_X^{\leq p-1})\]
(see \cite{GerVor, Saf} for more details).

Since we are in a local setting, we assume that $X$ is the spectrum of a polynomial algebra. In other words, we assume that $\D_X^{\leq p-1}=k\langle x_1,...,x_n, \del_1,...,\del_n\rangle$ where the relations are given by the usual Weyl-algebra relations with the extra relations $\del_i^p=0$.

Conider the following elements in $\D_X^{\leq p-1}\otimes_k \D_X^{\leq p-1}$:
\[\bar{\del}_i:=\del_i\otimes 1+1\otimes \del_i, \quad\mbox{and}\quad \bar{x}_j:=1\otimes x_j\]
where both $i$ and $j$ run from 1 to $n$. These elements, satisfy the usual Weyl-algebra relations
\[\bar{\del}_i\bar{x}_j-\bar{x}_j\bar{\del}_i=\delta_{ij}\]
and the relations $\bar{\del}_i^p=0$.

Using these elements, there is a unique algebra map $\rho: \D_X^{\leq p-1}\rightarrow \D_X^{\leq p-1}\otimes_k \D_X^{\leq p-1}$ sending the $\del_i$ to the $\bar{\del}_i$ and the $x_j$ to the $\bar{x}_j$. Under the projection 
\[\pi: \D_X^{\leq p-1}\otimes_k \D_X^{\leq p-1}\ra \D_X^{\leq p-1}\otimes_{\cO_X} \D_X^{\leq p-1}\]
these elements become the corresponding elements under the map $f: D_X\ra \D_X^{\leq p-1}\otimes_{\cO_X} \D_X^{\leq p-1}$ meaing that $\pi(\bar{\del}_i)=\phi(\del_i)$ and $\pi(\bar{x}_j)=\phi(x_j)$. This is true for any general element $\del \in \D_X^{\leq p-1}$, meaning that $(\pi\circ \rho)(\del)=\phi(\del)$. Indeed, let $\del=f\prod_{i\in I} \del_i$, we have
\[(\pi\circ\rho)(\del)=\sum_{J\subset I}\prod_{j\in J}\del_j \otimes f\prod_{i\in I\setminus J}\del_i.\]
Similarly, using the identification
\[\D_X^{\leq p-1}\otimes_{\cO_X}\D_X^{\leq p-1}\ra \sDiff_{\cO_{X'}}(\cO_X, \D_X^{\leq p-1})\]
we have that
\[\phi(\del)=(x\mapsto f\prod_{i\in I} \del_i \cdot x)=\sum_{J\subset I} (\prod_{j\in J} \del_j).x\cdot \prod_{i\in I\setminus J}\del_i\] showing that $(\pi\circ \rho)(\del)=\phi(\del)$.

Since  $\D_X^{\leq p-1}\iso \sDiff_{\cO_{X'}}(\cO_X^0,\D_X^{\leq p-1})$ as algebras, we obtain the corresponding elements $\bar{\del}_i$ and $\bar{x}_j$ inside $\sDiff_{\cO_{X'}}(\cO_X,\D_X^{\leq p-1}\otimes_{\cO_X}\D_X^{\leq p-1})$. Thus, locally, $\chi$ is a dg-algebra map, showing that globally it is is a dg-algebra map as well!

Finally, we show that $\chi$ is a quasi-isomorphism. First, we note that $Diff(\cO_X^\sbt, \D_X^{\leq p-1})^{op}$ has a natural filtration given by the filtration of $\D_X^{\leq p-1}$. This filtration is compatible with the $\sDiff_{\cO_{X'}}(\cO_X^\sbt, \cO_X)$-action. This induces a filtration on the derived tensor product 
\[\sDiff_{\cO_{X'}}(\cO_X^\sbt,\D_X^{\leq p-1, op})^{op}\otimes_{\sDiff_{\cO_{X'}}(\cO_X^\sbt,\cO_X)} \sDiff_{\cO_{X
}}(\cO_X^\sbt,\D_X^{\leq p-1}).\]
However, by Proposition \ref{prop:quant2}, we know that $\sDiff_{\cO_{X'}}(\cO_X^\sbt,(\D_X^{\leq p-1})^{op})^{op}$ is quasi-isomorphic to $\cO_X$, in other words, the graded pieces of the fitration are all concentrated in degree 0. As a consequence, the cohomology of the complex 
\[\cO_X \otimes^{\L}_{\sDiff_{\cO_{X'}}(\cO_X^\sbt,\cO_X)} \cO_X\]
has to be concentrated in degree 0 as well. We are done.\qed

\end{Proof}

\paragraph Consider the functors $F$ and $G$ between the coherent category of (left)-modules over $\cO_{X'}$ and over $\sDiff_{\cO_{X'}}(\cO_X^\sbt,\cO_X)$ given by:
\[F:\Coh(\sDiff_{\cO_{X'}}(\cO_X^\sbt,\cO_X))\ra \Coh(\cO_{X'})\quad F(-)=\cO_X\otimes_{\sDiff_{\cO_{X'}}(\cO_X^\sbt,\cO_X)}-\]
\[G:\Coh(\cO_{X'})\ra \Coh(\sDiff_{\cO_{X'}}(\cO_X^\sbt,\cO_X))\quad G(-)=\cO_X\otimes_{\cO_{X'}}.-\]
It is clear that the composite functor $FG(-)$ is given by $\D_X^{\leq p-1}\otimes_{\cO_{X'}}-$, hence for any algebra object $A$ in the coherent category $\Coh(\cO_{X'})$ we have that $FG(A)$ and $A$ are Morita equivalent. As a consequence for any object $B$ of the coherent category of left-modules over $\sDiff_{\cO_{X'}}(\cO_X^\sbt,\cO_X)$ which is of the form $G(A)$ we have $GF(B)$ is Morita-equivalent to $B$. This implies that $F$ and $G$ are essentially surjective functors between the $(\infty,2)$-categories of the coherent category of (left)-modules over $\cO_{X'}$ (where we consider $\cO_{X'}$ as an $E_2$-algebra) and of the full thick subcategory of the coherent category of (left)-modules over $\sDiff_{\cO_{X'}}(\cO_X,\cO_X)$ generated by $\cO_X$. This concludes the proof the Theorem \ref{thm:main1}.

\end{document}